\documentclass{article}
\usepackage{subeqn}

\usepackage{amsfonts}
\usepackage{graphicx}

\newtheorem{theorem}{Theorem}[section]
\newtheorem{example}{Example}[section]
\newtheorem{remark}{Remark}[section]
\newcommand{\dd}{{\, {\mathrm d}}}

\title{An Extension of the Well-Posedness Concept\\
for Fractional Differential Equations
of Caputo's Type}

\author{Kai Diethelm\thanks{Email: k.diethelm@tu-bs.de}
\thanks{\em Institut Computational Mathematics, Technische Universit\"at Braunschweig,
               Fallersleber-Tor-Wall 23, 38100 Braunschweig, Germany}
\thanks{\em GNS Gesellschaft f\"ur numerische Simulation mbH, Am Gau\ss berg 2,
               38114 Braunschweig, Germany}
 }

\begin{document}

\maketitle

\begin{abstract}
It is well known that, under standard assumptions, initial value
problems for fractional ordinary differential equations involving
Caputo-type derivatives are well posed in the sense that a unique
solution exists and that this solution continuously depends on the 
given function, the initial value and the order of the derivative.
Here we extend this well-posedness concept to the extent that we also
allow the location of the starting point of the differential operator
to be changed, and we prove that the solution depends on this
parameter in a continuous way too if the usual assumptions are
satisfied. 
Similarly, the solution to the corresponding terminal
value problems depends on the location of the starting point and
of the terminal point in a continuous way too.
\end{abstract}


Keywords:
Fractional differential equation, Caputo derivative, 
well-posed problem, continuous dependence, initial value,
terminal value, starting point

MSC 2010:
  34A08, 34A12, 26A33

\section{Introduction}

The goal of this brief note is to discuss a generalization of a
well known fundamental result from the theory of Caputo-type
fractional differential equations. Specifically, we are interested in the
classical initial value problem 
\begin{subequations}
  \label{eq:ivp}
  \begin{eqnarray}
    D_{*a}^\alpha y(t) &=& f(t, y(t)), 
    \label{eq:ivp-de} \\
    y^{(k)}(a) &=& y_k \qquad (k = 0, 1, \ldots, \lceil \alpha \rceil - 1)
    \label{eq:ivp-ic}
  \end{eqnarray}
\end{subequations}
for some $\alpha > 0$ on the interval $[a, T]$ 
where $D_{*a}^\alpha$ denotes the Caputo differential operator of
order $\alpha$ with starting point $a$ \cite[Chapter 3]{Diethelm2010}.
For this initial value problem, there holds the following existence and
uniqueness result \cite{DiethelmFord2002,Tisdell2012}: 

\begin{theorem}
  \label{thm:picard}
  Let $f:[a,T] \times \mathbb R \to \mathbb R$ be continuous and
  bounded, and assume that it satisfies a Lipschitz condition with
  respect to the second variable. Then, the initial value problem
  (\ref{eq:ivp}) has a unique continuous solution on $[a,T]$.
\end{theorem}

We also know that the solution $y$ of the initial value
problem (\ref{eq:ivp}) depends continuously on the given data
$f$, $y_k$ ($k = 0,1,\ldots, \lceil \alpha \rceil - 1$) and $\alpha$,
i.e.\ a small change in any of these values also implies a
small change in the solution 
(cf.\ \cite{DiethelmFord2002} or \cite[\S 6.3]{Diethelm2010}):

\begin{theorem}
  \label{thm:wellposed1}
  Let $y$ be the solution of (\ref{eq:ivp}), and let $\tilde y$ be the
  solution of the initial value problem 
  \begin{subequations}
    \label{eq:ivp2}
    \begin{eqnarray}
      D_{*a}^{\tilde \alpha} \tilde y(t) &=& \tilde f(t, \tilde y(t)), \\
      \tilde y^{(k)}(a) &=& \tilde y_k \qquad (k = 0, 1, \ldots, \lceil \tilde \alpha \rceil - 1)
    \end{eqnarray}
  \end{subequations}  
  where $f$ and $\tilde f$ are both assumed to satisfy the hypotheses
  of Theorem \ref{thm:picard}. Moreover, let $\lceil \alpha \rceil =
  \lceil \tilde \alpha \rceil$. Then, both initial value problems have
  unique continuous solutions $y$ and $\tilde y$, respectively, on
  $[a,T]$, and we have that  
  \begin{equation}
    \label{eq:wellposed1}
    \| y - \tilde y \|_\infty = O(\alpha - \tilde \alpha) + 
                                O(\| f - \tilde f\|_\infty) +
                                O(\max_k |y_k - \tilde y_k|).
  \end{equation}
\end{theorem}

In accordance with the terminology used in the classical case of
integer order differential equations, Theorem \ref{thm:wellposed1} is
usually summarized by saying that the initial value problem
(\ref{eq:ivp}) is well posed. This result is highly significant in
practical applications because it allows to conclude that a
mathematical model of the form (\ref{eq:ivp}) can provide useful
results --- i.e., results that differ from the correct values only by
a small amount --- even if the parameters of the process that is being
modeled are only known up to some limited accuracy.

However, there is a certain limitation in this theory because it only
allows to deal with small perturbations in the parameters appearing in
the function $f$, in the order of the differential operator $\alpha$
(together, these data typically describe material parameters or
similar properties), and in the initial values $y_k$ that describe
the state of the system at the start of the process. On the other
hand, the theory does not admit to investigate the behaviour of the
solution in the case of a small change of the starting point $a$ of
the differential operator in eq.\ (\ref{eq:ivp-de}), i.e.\ the point
at which the initial conditions (\ref{eq:ivp-ic}) are prescribed. It
is the goal of this paper to demonstrate that a small change in this
value also only leads to a small change in the solution, not only in
the neighbourhood of the starting point but throughout the complete
interval $[a,T]$ where the solution exists. Such a property is of
minor significance for modeling processes in a laboratory environment
where the starting time of the experiment is exactly known, but it can
be of utmost importance when mathematically simulating processes
observed in the real world where the starting time of the process is
known only unprecisely. Applications of this latter class include
many phenomena that have been successfully modeled using fractional
order equations such as, e.g. earthquakes
\cite{Caputo1967,CaputoMainardi1971,EvangelatosSpanos2011} (where
usually a good approximation, but not the exact value, of the starting
time is known) or the spreading of epidemics \cite{Diethelm2013} and
the distribution of pollutants in ground water
\cite{Benson2001,Meerschaert2011} where often at most a very rough
idea of the starting point exists.

\section{Main Result}
\label{sec:ivp}

The main result of our work is the following statement that allows us
to conclude that the solution to the initial value problem
(\ref{eq:ivp}) is indeed continuous with respect to the location of
the starting point.

\begin{theorem}
  \label{thm:main}
  Let $a \le \tilde a < T$, and consider the initial value problems
  (\ref{eq:ivp}) and
  \begin{subequations}
    \label{eq:ivp3}
    \begin{eqnarray}
      D_{*\tilde a}^\alpha \tilde y(t) &=& f(t, \tilde y(t)), 
      \label{eq:ivp3-de} \\
      \tilde y^{(k)}(\tilde a) &=& y_k \qquad (k = 0, 1, \ldots, \lceil \alpha \rceil - 1)
      \label{eq:ivp3-ic}
    \end{eqnarray}
  \end{subequations}
  under the assumptions of Theorem \ref{thm:picard}. 
  The solutions $y$ and $\tilde y$ to these
  initial value problems satisfy the relation
  \begin{equation}
    \label{eq:wellposed2}
    \sup_{t \in [\tilde a, T]} | y(t) - \tilde y (t)| = O(|a - \tilde a|^{\min \{ \alpha, 1\}}) .
  \end{equation}
\end{theorem}

\begin{remark}
  The only difference between problems (\ref{eq:ivp}) and
  (\ref{eq:ivp3}) is the location of the starting point of the
  fractional differential operator in the differential equation that,
  since we are talking about an initial value problem, coincides with
  the point at which the initial condition is prescribed. This is
  sufficient because the effects of perturbations in all other
  parameters are already known from Theorem \ref{thm:wellposed1}.
\end{remark}

\begin{remark}
  As one can, if necessary, always exchange the roles of $y$ and $\tilde y$,
  the assumption that $a \le \tilde a$ that we have imposed in Theorem
  \ref{thm:main} does not imply a loss of generality.
\end{remark}

\begin{remark}
  By construction, the solution $y$ of the problem (\ref{eq:ivp}) is
  defined on the interval $[a,T]$. Similarly, the function $\tilde y$
  that solves the initial value problem (\ref{eq:ivp3}) is defined on
  $[\tilde a, T]$. It is therefore perfectly natural to perform the
  comparison of the two functions $y$ and $\tilde y$ in eq.\
  (\ref{eq:wellposed2}) on the intersection of these two intervals,
  i.e.\ (in view of the assumption $a \le \tilde a$ that we
  had imposed in Theorem \ref{thm:main}) on the interval
  $[\tilde a, T]$.
\end{remark}

{\bf Proof.}
It is well known \cite{DiethelmFord2002} that the initial value
problem (\ref{eq:ivp}) is equivalent to the Volterra integral equation
\begin{equation}
  \label{eq:volterra1}
  y(t) = \sum_{k=0}^{\lceil \alpha \rceil - 1} \frac{y_k}{k!} (t - a)^k
          + \frac1{\Gamma(\alpha)} \int_a^t (t-s)^{\alpha-1} f(s, y(s)) \dd s
\end{equation}
for $t \in [a, T]$, and similarly (\ref{eq:ivp3}) is equivalent to 
\begin{equation}
  \label{eq:volterra2}
  \tilde y(t) = \sum_{k=0}^{\lceil \alpha \rceil - 1} \frac{y_k}{k!} (t - \tilde a)^k
          + \frac1{\Gamma(\alpha)} \int_{\tilde a}^t (t-s)^{\alpha-1} f(s, \tilde y(s)) \dd s
\end{equation}
for $t \in [\tilde a, T] \subseteq [a,T]$. Thus, subtracting
(\ref{eq:volterra2}) from (\ref{eq:volterra1}), we obtain for $t \in
[\tilde a, T]$
\begin{subequations}
  \begin{equation}
    \label{eq:diff1}
    \delta(t) := |y(t) - \tilde y(t)| \le |D_1| + |D_2| + |D_3|
  \end{equation}
  where
  \begin{eqnarray}
    \label{eq:diff21}
    D_1 &=& \sum_{k=1}^{\lceil \alpha \rceil - 1} \frac{y_k}{k!} [(t - a)^k - (t - \tilde a)^k], \\
    \label{eq:diff22}
    D_2 &=&  \frac1{\Gamma(\alpha)} \int_a^{\tilde a} (t-s)^{\alpha-1} f(s, y(s)) \dd s, \\
    \label{eq:diff23}
    D_3 &=& \frac1{\Gamma(\alpha)} \int_{\tilde a}^t (t-s)^{\alpha-1} [f(s, y(s)) - f(s, \tilde y(s))] \dd s.
  \end{eqnarray}
\end{subequations}

The Lipschitz condition on $f$ then implies, if we denote the
corresponding Lip\-schitz constant by $L$, that
\begin{equation}
  \label{eq:diff33}
  |D_3| \le \frac L{\Gamma(\alpha)} \int_{\tilde a}^t (t-s)^{\alpha-1} \delta(s) \dd s.
\end{equation}

We now need to distinguish two cases. In the first case, $0 < \alpha
\le 1$, we clearly have
\[
  D_1 = 0
\]
and, denoting by $M$ the supremum of $f$ on its domain of definition
(which, by assumption, is finite),
\[
  |D_2| \le \frac M{\Gamma(\alpha)} \int_a^{\tilde a} (t - s)^{\alpha-1} \dd s 
        \le \frac M{\Gamma(\alpha)} \int_a^{\tilde a} (\tilde a - s)^{\alpha-1} \dd s 
	 =  \frac {M(\tilde a - a)^\alpha}{\Gamma(\alpha+1)} .
\]

In the other case, $\alpha > 1$, the mean value theorem of
differential calculus implies
\[
  |D_1| \le |\tilde a - a| \sum_{k=1}^{\lceil \alpha \rceil - 1} \frac{y_k}{(k-1)!} (t - \xi_k)^{k-1}
\]
with certain $\xi_k \in [a, \tilde a]$, and hence
\[
  |D_1| \le |\tilde a - a| \sum_{k=1}^{\lceil \alpha \rceil - 1} \frac{y_k}{(k-1)!} (T-a)^{k-1}
         =  C |\tilde a - a|
\]
where $C$ is a constant independent of $t$ and $a - \tilde
a$. Moreover, in this case we have, with $M$ again denoting the
supremum of $f$, that
\begin{eqnarray*}
  |D_2| & \le & \frac M{\Gamma(\alpha)} \int_a^{\tilde a} (t - s)^{\alpha-1} \dd s 
          \le \frac M{\Gamma(\alpha)} \int_a^{\tilde a} (T - a)^{\alpha-1} \dd s \\
	&  =  & \frac {M}{\Gamma(\alpha)} (\tilde a - a) (T - a)^{\alpha-1} .
\end{eqnarray*}

It now follows from eqs.\ (\ref{eq:diff1}) and (\ref{eq:diff33}) that
\begin{equation}
  \label{eq:gronwall}
  \delta(t) \le |D_1| + |D_2| + \frac L{\Gamma(\alpha)} \int_{\tilde a}^t (t-s)^{\alpha-1} \delta(s) \dd s,
\end{equation}
and our estimates above imply that, in either case, 
\[ 
  |D_1| + |D_2| \le C^* |a - \tilde a|^p \qquad \mbox{with } p := \min\{\alpha,1\}
\]
where $C^*$ is an absolute constant. Thus, denoting (as usual) the one-parameter
Mittag-Leffler function by $E_\alpha(z) := \sum_{k=0}^\infty z^k/\Gamma(\alpha k + 1)$, 
the fractional version of Gronwall's Lemma \cite[Lemma 6.19]{Diethelm2010} 
allows us to conclude from eq.\ (\ref{eq:gronwall}) that
\[
  | y(t) - \tilde y(t)|  =  \delta(t) \le C^* |a - \tilde a|^p E_\alpha (L (t-\tilde a)^\alpha) 
                        \le  C^* |a - \tilde a|^p E_\alpha (L (T-a)^\alpha)
\]
for all $t \in [\tilde a, T]$ which, since 
$C^* E_\alpha (L (T-a)^\alpha)$ is independent of $t$, completes the proof.
~$\Box$

We point out that the estimate of eq.\ (\ref{eq:wellposed2}) is best possible:

\begin{theorem}
  Under the assumptions of Theorem \ref{thm:main}, the estimate of eq.\
  (\ref{eq:wellposed2}) cannot be improved, i.e.\ the exponent $\min\{\alpha,1\}$ in
  the $O$-term on the right-hand side of eq.\ (\ref{eq:wellposed2}) cannot be
  replaced by a larger number. 
\end{theorem}

In order to prove this statement, it suffices to construct examples where the order
indicated in this $O$-term is actually attained. We shall do this separately
for the two cases $0 < \alpha \le 1$ and $\alpha > 1$:

\begin{example}
  For the case $0 < \alpha \le 1$, we set $a = 0$ and consider the initial
  value problems 
  \begin{equation}
    \label{eq:ex11}
    D_{*0}^\alpha y(t) = y(t), \qquad y(0) = 1,
  \end{equation}
  and
  \begin{equation}
    \label{eq:ex12}
    D_{*\tilde a}^\alpha \tilde y(t) = \tilde y(t), \qquad \tilde y(\tilde a) = 1
  \end{equation}
  whose solutions are well known to be 
  \begin{equation}
    \label{eq:ex1s}
    y(t) = E_\alpha (t^\alpha)   \qquad\mbox{ and }\qquad 
    \tilde y(t) = E_\alpha ((t-\tilde a)^\alpha)
  \end{equation}
  where again $E_\alpha(z)$ denotes the one-parameter Mittag-Leffler function.
  Here we easily conclude (for arbitrary $T > \tilde a$)
  \begin{eqnarray*}
    \sup_{t \in [\tilde a, T]} |y(t) - \tilde y(t)|
    & \ge & y(\tilde a) - \tilde y(\tilde a) 
       =    E_\alpha ({\tilde a}^\alpha) - 1
       =    \sum_{k=1}^\infty \frac{{\tilde a}^{\alpha k}}{\Gamma(\alpha k + 1)} \\
    & \ge & \frac{{\tilde a}^{\alpha}}{\Gamma(\alpha + 1)}
       =    \frac{|a - \tilde a|^{\alpha}}{\Gamma(\alpha + 1)}
  \end{eqnarray*}
  which, since in this case $\alpha = \min\{\alpha,1\}$, gives the required
  result. 
  \label{expl1}
\end{example}

\begin{example}
  For the case $\alpha > 1$, we use the same initial value problems as 
  in Example \ref{expl1} above and choose an arbitrary $T > \tilde a$. Then,
  using the mean value theorem, we deduce
  \begin{eqnarray*}
    \sup_{t \in [\tilde a, T]} |y(t) - \tilde y(t)|
    & \ge & y(T) - \tilde y(T) 
       =    E_\alpha (T^\alpha) - E_\alpha ((T-\tilde a)^\alpha) \\
    &  =  & \sum_{k=1}^\infty \frac{T^{\alpha k} - (T-\tilde a)^{\alpha k}}{\Gamma(\alpha k + 1)} 
      \ge   \frac{T^{\alpha} - (T-\tilde a)^{\alpha}}{\Gamma(\alpha + 1)}  
       =    \frac{\xi^{\alpha-1}}{\Gamma(\alpha)} \tilde a \\
    & \ge & \frac{(T-\tilde a)^{\alpha-1}}{\Gamma(\alpha)} | a - \tilde a|
  \end{eqnarray*}
  which, since in this case $1 = \min\{\alpha,1\}$, again gives the required
  result.   
\end{example}
  
\begin{remark}
  The investigation of questions of this type is relevant in
  connection with fractional differential equations mainly because the
  associated operators exhibit a certain memory. It is also possible
  to develop other memory-dependent operators, and for the
  corresponding operator equations one would then need to look at the
  same type of questions. We believe that the technique
  employed here will be useful in those settings too.
\end{remark}

\section{Terminal Value Problems}

So far, we have discussed the question of the dependence of the
solution to a fractional-order initial value problem on the location of the
starting point, and we have seen that, under reasonable assumptions, this
dependence is of a continuous nature. A related question is whether the same
type of dependence can be proved for terminal value problems. In this case one
would first investigate whether the solution $y$ to the problem
$$
  D_{*a}^\alpha y(t) = f(t, y(t)), \qquad
  y^{(k)}(T) = y_k \quad (k = 0, 1, \ldots, \lceil \alpha \rceil - 1)
$$ 
where $T > a$ depends on $a$, and possibly also on $T$, in a continuous way. 

The existence and uniqueness of solutions to such problems and some related
questions have been addressed in \cite{Diethelm2008,DiethelmFord2012}, and is
has turned out that the case $0 < \alpha < 1$ is of particular interest. 
We shall therefore concentrate on this case.  The question for
the well-posedness of such problems in the classical sense, i.e.\ if $f$, $\alpha$
and the initial values $y_k$ are varied, has been addressed in \cite{FordMorgado2013}.
In the context of the problem under consideration here, it is then very natural
to consider $a$ as an additional unknown and to ask under which conditions it
is possible to identify some suitable additional information given which one
can conclude that it is possible to uniquely determine both the solution $y$
to the terminal value problem and the starting point $a$.  Another
topic of interest in such a connection would be to find out how the solution
reacts to a small change of the value $T$.

Questions like these arise rather naturally if one tries to model a phenomenon
observed outside of a laboratory by a fractional differential equation and
does not know when the process has started. An example could be, e.g., the
case where the dynamics of an epidemic are discussed as in \cite{Diethelm2013}.
In such a case it is quite natural to have observations describing the
development of the disease at certain points of time, but typically not at the
instant where the infection first reached the population and started the
process because this point in time is simply not known. Then, one is typically
interested in finding out how the epidemic progresses, and in order to answer
this question the model requires that one first determines the starting point.

While the extension of our results above to this area might 
at first seem to be a straightforward matter, a closer look reveals 
certain significant differences. Most notably, while the initial value 
problems discussed in Section \ref{sec:ivp} were equivalent to integral 
equations of Volterra's type, cf.\ eq.\ (\ref{eq:volterra1}), the natural 
integral equation formulation
of the terminal value problems under consideration now has a Fredholm form 
\cite[Theorem 6.18]{Diethelm2010}. Nevertheless it is possible to show 
similar results in this case too. We shall first state the result for the 
case that $T$ varies.

\begin{theorem}
  \label{thm:tvp1}
  Consider the terminal value problems
  \begin{equation}
    \label{eq:tvp1}
    D_{*a}^\alpha y(t) = f(t, y(t)), \qquad
    y(T) = y^*,
  \end{equation}
  and
  \begin{equation}
    \label{eq:tvp2}
    D_{*a}^\alpha \tilde y(t) = f(t, \tilde y(t)), \qquad
    \tilde y(\tilde T) = y^*,
  \end{equation}
  for some $\alpha \in (0,1)$ and $a < T \le \tilde T$,
  where the given function $f$ is once again
  assumed to satisfy the conditions of Theorem \ref{thm:picard}.
  The solutions $y$ and $\tilde y$ to these terminal value problems
  satisfy the relation
  \begin{equation}
    \label{eq:tvp-bound1}
    \sup_{t \in [a, T]} | y(t) - \tilde y(t) | = O(|T - \tilde T|^\alpha).
  \end{equation}
\end{theorem}

{\bf Proof.}
  Following \cite[Theorem 6.18]{Diethelm2010}, we rewrite the terminal value problems
  (\ref{eq:tvp1}) and (\ref{eq:tvp2}) as equivalent integral equations,
  \begin{equation}
    \label{eq:fredholm1}
    y(t) = y^* + \frac1{\Gamma(\alpha)} \int_a^T G(t,s) f(s, y(s)) \dd s
  \end{equation}
  and
  \begin{equation}
    \label{eq:fredholm2}
    \tilde y(t) = y^* + \frac1{\Gamma(\alpha)} \int_a^{\tilde T} \tilde G(t,s) f(s, \tilde y(s)) \dd s,
  \end{equation}
  respectively, where
  $$
    G(t,s) = 
    \cases{
      - (T-s)^{\alpha-1}                  & \mbox{for $s > t$}, \cr
      (t-s)^{\alpha-1} - (T-s)^{\alpha-1} & \mbox{for $s \le t$}, \cr
    }
  $$
  and an analog relation, with $T$ being replaced by $\tilde T$, holds for $\tilde G$.
  Subtracting (\ref{eq:fredholm2}) from (\ref{eq:fredholm1}) and recalling that $a < T \le \tilde T$, 
  we obtain
  \begin{eqnarray} 
    \label{eq:tvp-proof1a}
    y(t) - \tilde y(t) 
    &=& \frac1{\Gamma(\alpha)} \int_a^T [G(t,s) f(s, y(s)) - \tilde G(t,s) f(s, \tilde y(s))] \dd s \\
    \nonumber
    & & {}   - \frac1{\Gamma(\alpha)} \int_T^{\tilde T} \tilde G(t,s) f(s, \tilde y(s)) \dd s
  \end{eqnarray}
  for $a \le t \le T$. In view of this inequality and the corresponding branch of the
  definition of $\tilde G$, the second integral in the representation (\ref{eq:tvp-proof1a}) 
  can be bounded in modulus as
  $$
    \left| \int_T^{\tilde T} \tilde G(t,s) f(s, \tilde y(s)) \dd s \right|
      \le \| f \|_\infty \int_T^{\tilde T} (\tilde T - s)^{\alpha-1} \dd s
       =  \frac{\| f \|_\infty}{\alpha} (\tilde T - T)^\alpha. 
  $$
  The remaining integral in eq.\ (\ref{eq:tvp-proof1a}) needs to be split up once again.
  This yields
  \begin{eqnarray*}
    \lefteqn{\int_a^T [G(t,s) f(s, y(s)) - \tilde G(t,s) f(s, \tilde y(s))] \dd s} \\
    & = & \int_a^T G(t,s) [f(s, y(s)) - f(s, \tilde y(s))] \dd s
          + \int_a^T [G(t,s) - \tilde G(t,s)] f(s, \tilde y(s)) \dd s.
  \end{eqnarray*}
  The second of these integrals can be bounded in modulus by
  \begin{eqnarray*}
    \lefteqn{\left | \int_a^T [G(t,s) - \tilde G(t,s)] f(s, \tilde y(s)) \dd s \right| } \\
    & \le & \| f \|_\infty \int_a^T |G(t,s) - \tilde G(t,s)| \dd s 
       =    \| f \|_\infty \int_a^T | (T-s)^{\alpha-1} - (\tilde T - s)^{\alpha -1}| \dd s \\
    &  =  & \frac1\alpha \| f \|_\infty [(T-a)^\alpha - (\tilde T - a)^\alpha + (\tilde T - T)^\alpha]
      \le   \frac1\alpha \| f \|_\infty (\tilde T - T)^\alpha.
  \end{eqnarray*}
  Combining the estimates obtained so far with the Lipschitz property of $f$ with respect
  to the second variable, we arrive at
  $$
    | y(t) - \tilde y(t) | 
      \le \frac{2 \| f \|_\infty}{\Gamma(\alpha+1)} (\tilde T - T)^\alpha
             + \frac L{\Gamma(\alpha)} \int_a^T |G(t,s)| \cdot |y(s) - \tilde y(s)| \dd s
  $$
  In order to conclude the desired inequality (\ref{eq:tvp-bound1}) from this relation,
  we need to invoke a Gronwall type argument. A suitable result of this type can be
  derived from \cite[Theorem 2.1]{Muresan1999} by noting that, as in the argumentation of 
  \cite[Theorems 3.13 and 4.8]{Linz1985}, the continuity requirement
  for the kernel function that is present in \cite{Muresan1999} can be relaxed to a 
  weaker integrability condition satisfied by our functions $G$ and $\tilde G$. 
~$\Box$

The case that $a$ varies can be handled in a similar (but not
exactly identical) way, and we can show the following result that is formally 
essentially the same as the previous theorem.

\begin{theorem}
  Consider the terminal value problems
  \begin{equation}
    \label{eq:tvp3}
    D_{*a}^\alpha y(t) = f(t, y(t)), \qquad
    y(T) = y^*,
  \end{equation}
  and
  \begin{equation}
    \label{eq:tvp4}
      D_{*\tilde a}^\alpha \tilde y(t) = f(t, \tilde y(t)), \qquad
      \tilde y(T) = y^*,
  \end{equation}
  for some $\alpha \in (0,1)$ and $a \le \tilde a < T$, 
  where the given function $f$ is once again 
  assumed to satisfy the conditions of Theorem \ref{thm:picard}.
  The solutions $y$ and $\tilde y$ to these terminal value problems
  satisfy the relation
  \begin{equation}
    \label{eq:tvp-bound2}
    \sup_{t \in [\tilde a, T]} | y(t) - \tilde y(t) | = O(|a - \tilde a|^\alpha).
  \end{equation}
\end{theorem}

{\bf Proof.}
  The basic steps are quite similar to the proof of Theorem \ref{thm:tvp1}. We begin
  by setting up the Fredholm equations for the terminal value problems (\ref{eq:tvp3})
  and (\ref{eq:tvp4}) and subtract them from each other. We note that the kernel functions
  appearing there depend only on $T$ and not on $a$. Therefore, as $T$ does not change 
  in the present setting, we have the same kernel $G$ in both integral equations; this 
  simplifies the analysis. Specifically, we obtain
  \begin{eqnarray}
    \label{eq:tvp-proof2a}
    |y(t) - \tilde y(t)|  & \le &  
      \frac1{\Gamma(\alpha)} \int_a^{\tilde a} |G(t,s) f(s, y(s))| \dd s \\
     \nonumber
     && {}
      + \frac1{\Gamma(\alpha)} \int_{\tilde a}^T |G(t,s)| \cdot |f(s, y(s)) - f(s, \tilde y(s))| \dd s 
  \end{eqnarray}
  The first integral in eq.\ (\ref{eq:tvp-proof2a}) can be estimated by
  \begin{eqnarray*}
    \lefteqn{\frac1{\Gamma(\alpha)} \int_a^{\tilde a} |G(t,s) f(s, y(s))| \dd s } \\
      & \le & \frac{\| f \|_\infty}{\Gamma(\alpha)} \int_a^{\tilde a} | G(t,s)| \dd s 
         =    \frac{\| f \|_\infty}{\Gamma(\alpha)} \int_a^{\tilde a} | (t-s)^{\alpha-1} - (T-s)^{\alpha-1}| \dd s \\
      &  =  & \frac{\| f \|_\infty}{\Gamma(\alpha+1)} [ (t-a)^\alpha - (t - \tilde a)^\alpha
                                                        - (T-a)^\alpha + (T-\tilde a)^\alpha] \\
      & \le & C_1 |a - \tilde a| + C_2 |a - \tilde a|^\alpha \le C |a - \tilde a|^\alpha
  \end{eqnarray*}
  with certain constants $C_1$, $C_2$ and $C$ because of the mean value theorem of differential
  calculus and the fact that the function $(\cdot)^\alpha$ satisfies a H\"older condition 
  of order $\alpha$.
  For the second integral in eq.\ (\ref{eq:tvp-proof2a}) we can simply use the Lipschitz property 
  of $f$ with respect to the second variable. Combining these two estimates for the two
  integrals in eq.\ (\ref{eq:tvp-proof2a}), we arrive at
  $$
    |y(t) - \tilde y(t)| 
       \le C |a - \tilde a|^\alpha 
           + \frac L {\Gamma(\alpha)} \int_{\tilde a}^T |G(t,s)| \cdot |y(s) - \tilde y(s)| \dd s.
  $$
  From here, we may again (as in the proof of Theorem \ref{thm:tvp1})
  argue with the help of the Gronwall inequality for Fredholm operators
  and obtain the result (\ref{eq:tvp-bound2}).
~$\Box$

We shall address additional questions related to problems 
of this sort, and the quest for numerical methods
for their solution (in this context, see \cite{FordMR} for first 
results), in a forthcoming separate paper.



\end{document}